\documentclass[10pt,leqno]{article}

\usepackage{amsmath,amssymb,amsthm,mathrsfs,mathtools}
\usepackage[margin=3cm]{geometry}
\usepackage{titlesec}
\usepackage{microtype}
\usepackage{indentfirst}
\usepackage{needspace}
\usepackage{fancyhdr}
\usepackage[hidelinks]{hyperref}
\titleformat{\subsection}{\normalfont\itshape}{\thesubsection.\enspace}{1pt}{}
\newtheorem{theorem}{Theorem}[section]
\newtheorem{lemma}[theorem]{Lemma}

\newtheorem{proposition}[theorem]{Proposition}
\theoremstyle{remark}
\newtheorem{remark}[theorem]{Remark}
\numberwithin{equation}{section}
\allowdisplaybreaks
\hypersetup{
  pdftitle={A one-component Liouville theorem for stationary magnetohydrodynamics in weak Lebesgue spaces},
  pdfauthor={Cuiyang Liu and Wei Luo},
  pdfkeywords={stationary MHD, Liouville theorem, one component, weak Lebesgue spaces, mixed norms, Hardy space}
}

\begin{document}
\title{A one-component Liouville theorem for \\ stationary MHD system in weak Lebesgue spaces}
\author{Cuiyang $\mbox{Liu}^1$\footnote{E-mail: liucy25@mail2.sysu.edu.cn},\quad
Wei $\mbox{Luo}^2$\footnote{E-mail: luow23@mail.sysu.edu.cn}\\
$^1\mbox{School}$ of Science, Shenzhen Campus of Sun Yat-sen University,\\
Shenzhen 518107, China\\
$^2\mbox{Department}$ of Mathematics, Sun Yat-sen University,\\
Guangzhou 510275, China}
\date{}
\maketitle
\hrule

\begin{abstract}
We prove a one-component Liouville theorem for stationary incompressible MHD system in $\mathbb R^3$: for smooth solutions with finite Dirichlet integral that decay at infinity, $u_3\in L^{p,\infty}$ and $b_3\in L^{q,\infty}$ with $1\le p,q<2$ force $u=b=0$ and constant pressure.
We also consider the case where the viscosity equals to the magnetic diffusivity, which needs fewer assumption since it has a simpler structure.
The core of our proof is a structural lemma adapted to MHD which replaces the head-pressure structure of the Navier--Stokes equations. It shows the reason why we can control the pressure with only one component.

\vspace*{5pt}

\vspace*{5pt}
\end{abstract}
\vspace*{10pt}
\begingroup
\linespread{1}\selectfont
\tableofcontents
\endgroup
\clearpage

% ===========================================================================
% Section 1 -- Introduction
% (formerly chapters/01_introduction_v2.tex)
% ===========================================================================
% Chapter 1 -- Introduction (Section 1 of the paper)
\section{Introduction}
\label{sec:introduction}

In this paper, We consider the stationary incompressible magnetohydrodynamic(MHD) system in $\mathbb{R}^3$
\begin{equation}\label{eq:mhd}
\begin{cases}
-\nu\Delta u+(u\cdot\nabla)u-(b\cdot\nabla)b+\nabla P=0,\\
-\eta\Delta b+(u\cdot\nabla)b-(b\cdot\nabla)u=0,\\
\operatorname{div}u=\operatorname{div}b=0,
\end{cases}
\end{equation}
where $u$ is the velocity, $b$ is the magnetic field, and $P$ is the total pressure, including the mechanical part and the magnetic part.
The constants $\nu>0$ and $\eta>0$ denote the viscosity and the magnetic diffusivity, which are not need to be equal.
We write $x=(x_h,z)$ with $x_h=(x_1,x_2)$, and $u_h=(u_1,u_2)$, $b_h=(b_1,b_2)$.
Throughout this paper, the solution is smooth and satisfies
\begin{equation}\label{eq:background}
\int_{\mathbb R^3}(|\nabla u|^2+|\nabla b|^2)\,dx<\infty,
\qquad
\lim_{R\to\infty}\sup_{|x|>R}(|u(x)|+|b(x)|)=0.
\end{equation}
This means that the solution has finite energy, and vanishes at infinity.

Our main result is to show that a condition controlling just one component of each field already forces both fields to be identically zero. In detail, if the third components of both fields lie in weak Lebesgue spaces of exponent below $2$ strictly, then both $u$ and $b$ vanishes.
Here the weak Lebesgue spaces is defined by the following quasi-norm, for $1\le p<\infty$ we write
\[
\|v\|_{L^{p,\infty}}=\sup_{t>0}t\,\bigl|\{x\in\mathbb R^3:|v(x)|>t\}\bigr|^{1/p}.
\]
The precise statements are Theorems~\ref{thm:main} and~\ref{thm:elsasser-single} below.

For suffiently regular (e.g. compactly supported smooth) fields, the energy identity can be computed formally.
Testing the two equations by $u$ and $b$ and summing them, the convective terms vanish because the free-divergence conditions, and the mixed terms cancel in pairs that is:
\[
\int_{\mathbb R^3}\bigl[(b\cdot\nabla)b\cdot u+(b\cdot\nabla)u\cdot b\bigr]\,dx
=\int_{\mathbb R^3}b\cdot\nabla(u\cdot b)\,dx=0.
\]
Hence, we obtain the equality $\nu\|\nabla u\|_2^2+\eta\|\nabla b\|_2^2=0$.
This implies $u=b=0$.
However, under \eqref{eq:background} this computation cannot be applied on $\mathbb R^3$: we need to use cutoff functions unavoidably, which produce new integral terms. Whether these new terms tend to zero is a nontrivial problem.

By now, a number of Liouville-type theorems have been established for the MHD system and its related Hall--MHD system.
For example, Chen, Li and Wang~\cite{ChenLiWang} obtained a criterion with $u\in\mathrm{BMO}^{-1}$ and $b\in L^{6,\infty}$; and Cho, Neustupa and Yang~\cite{ChoNeustupaYang2024} got further criteria on Lebesgue norms in form of growth conditions; both papers also consider the Hall--MHD system.
Chae~\cite{Chae2023} proved anisotropic Liouville theorems for MHD of all dimensions.

The classic stationary Navier--Stokes equations are the comparison problem: when $b=0$, \eqref{eq:mhd} reduces to the stationary Navier--Stokes system in $\mathbb{R}^3$, for which the same question remains open for a long time; see Galdi~\cite{Galdi}.
For the velocity, the Galdi's condition $u\in L^{9/2}(\mathbb R^3)$ implies $u=0$, which can be viewed as the one of earliest results on this topic; see \cite[Theorem X.9.5]{Galdi}.
Chae and Wolf~\cite{ChaeWolf2016} relaxed it in a logarithmical way, Seregin~\cite{Seregin2016} used a novel $u\in L^6\cap\mathrm{BMO}^{-1}$ condition, while Jarr\'in~\cite{Jarrin2020} developed Lorentz and Morrey criteria with a broader range.
For the $L^{9/2, \infty}$ assumption, which can be seen as a direct weaken of Galdi's result, Kozono, Terasawa and Wakasugi~\cite{KTW2017} need some smallness conditions relative to the Dirichlet integral, but whether the pure weak $L^{9/2}$ condition suffices is still open.

Another main line of work is using anisotropic or componentwise assumption.
For the example of anisotrpic criteria cases, we can see Luo, Yin~\cite{LuoYin2017}; Phan~\cite{Phan}, Chae~\cite{ChaeNS2023} and Zhang and Zu~\cite{ZhangZu2025}.
Meanwhile, Carrillo, Pan, Zhang and Zhao~\cite{CarrilloPanZhangZhao2020} and Tsai~\cite{Tsai2021} showed how symmetry, periodicity, boundary conditions, or annular control affect the energy.
Vergara-Hermosilla~\cite{VergaraHermosilla2026} obtained a Liouville theorem from the radial velocity component.
In our earlier paper~\cite{LiuLuoNS} we proved vanishing under $u_3\in L^{3/2,\infty}$ for the Navier--Stokes system.

This paper gives a sufficient condition of this componentwise type for the MHD system, for arbitrary $\nu,\eta>0$.
Only the third component of each field is controlled additionally, the other two exponents are independent.
The endpoint is the exponent $2$: whether $2$ is optimal is left open, and Remark~\ref{rem:endpoint} discusses it detailedly.
In the case of equal diffusivity i.e. $\nu = \eta$, we use the Els\"asser formulation yield a single component of one Els\"asser variable is enough.

\subsection{Main results}

All spaces in the sequel are over $\mathbb R^3$ unless a different domain is indicated.

\begin{theorem}\label{thm:main}
Let $(u,b,P)$ be a smooth solution of \eqref{eq:mhd} satisfying \eqref{eq:background}.
Suppose that, for some $p,q\in[1,2)$,
\begin{equation}\label{eq:main-assumption}
u_3\in L^{p,\infty},\qquad b_3\in L^{q,\infty}.
\end{equation}
Then $u=b=0$ and $P$ is constant.
\end{theorem}

In the assumptions of above theorem, the exponents can be completely independent, and no any smallness is needed.
The endpoints $p=2$ and $q=2$ are not covered; see Remark~\ref{rem:endpoint}.
When $b=0$, Theorem~\ref{thm:main} gives a Navier--Stokes criterion under $u_3\in L^{p,\infty}$ for every $1\le p<2$, which contains the exponent $3/2$ of~\cite{LiuLuoNS}.

Set $\zeta^\pm=u\pm b$.
When the viscosity and the magnetic diffusivity coincide i.e. $\nu = \eta$,  one component of a single Els\"asser field has been sufficient.

\begin{theorem}\label{thm:elsasser-single}
Let $(u,b,P)$ be a smooth solution of \eqref{eq:mhd} satisfying \eqref{eq:background}, and suppose that $\nu=\eta=\kappa>0$.
If, for some $p\in[1,2)$, either
\begin{equation}\label{eq:elsasser-single-assumption}
\zeta_3^+=u_3+b_3\in L^{p,\infty}
\qquad\text{or}\qquad
\zeta_3^-=u_3-b_3\in L^{p,\infty},
\end{equation}
then $u=b=0$ and $P$ is constant.
\end{theorem}

In the equal-diffusivity case this hypothesis is weaker than that of Theorem~\ref{thm:main}. This is natural, since this case is intermediate between Navier-Stokes system and our \eqref{eq:mhd}.
The restriction $\nu=\eta$ removes the cross-diffusion term, and the endpoint $p=2$ is again not claimed.

\subsection{Strategy of the proof}

For Navier--Stokes system, a usual idea is using the head-pressure : $P$ is paired with $\frac12|u|^2$, so that this quantity satisfies a second order elliptic equation then satisfies a maximum principle, from which its sign can be determined.
A direct but complicated computation can show that no quantity of the form $p+Q(|u|,|b|)$,
with $Q$ a quadratic polynomial of $|u|$ and $|b|$,
can satisfy a second order elliptic equation.
This fact means that we can not use the method in Navier-Stoke system to get our MHD results directly. 
So, we need to find another structure.
It has the following form.

\begin{lemma}\label{lem:structural}
Let $(u,b,P)$ solve \eqref{eq:mhd} and define
\[
f=\nu u_3,\qquad F=-(u_3u_h-b_3b_h),\qquad
E=\nu\partial_zu_3-u_3^2+b_3^2.
\]
Then
\begin{equation}\label{eq:intro-directional}
\partial_z(P-E)=\Delta_h f+\operatorname{div}_h F.
\end{equation}
\end{lemma}

\begin{proof}
The third component of \eqref{eq:mhd} reads
$-\nu\Delta u_3+(u\cdot\nabla)u_3-(b\cdot\nabla)b_3+\partial_zP=0$.
Because $\operatorname{div}u=\operatorname{div}b=0$,
\[
(u\cdot\nabla)u_3=\partial_z(u_3^2)+\operatorname{div}_h(u_3u_h),
\qquad
(b\cdot\nabla)b_3=\partial_z(b_3^2)+\operatorname{div}_h(b_3b_h).
\]
Substituting these and splitting $\Delta=\Delta_h+\partial_z^2$, the terms carrying a vertical derivative collect on the left:
\[
\partial_z\bigl(P-\nu\partial_zu_3+u_3^2-b_3^2\bigr)
=\Delta_h(\nu u_3)-\operatorname{div}_h(u_3u_h-b_3b_h),
\]
which is \eqref{eq:intro-directional}.
\end{proof}

In the above equation, we subtract $E$ from $P$. This is because $E$ is built from the third components $u_3$ and $b_3$,
which are precisely the quantities we give extra assumptions.
Hence, $E$ does not require structural treatment.

Aftering building the \eqref{eq:intro-directional}, we will use a method of heat-kernel decomposition
into high and low frequencies, combined with an iteration, to obtain the
structural estimate Proposition~\ref{prop:pressure} of this identity.
Then, it remains to Check term by term that the quantities in \eqref{eq:intro-directional} satisfy
the conditions of that estimate. Then yields the Liouville theorem.
The argument uses the div--curl theorem of Coifman, Lions, Meyer and Semmes and
the Hardy-space theory of compensated compactness quantities~\cite{CLMS,Lindberg}.

Both theorems are proved by applying Proposition~\ref{prop:pressure} to different data, $(u_3,b_3)$ in the first case and $(\zeta_3^+,\zeta_3^-)$ in the second.

% ===========================================================================
% Section 2 -- Preliminaries
% (formerly chapters/02_preliminaries.tex)
% ===========================================================================
\section{Preliminaries}\label{sec:prelim}

\subsection{Weak interpolation, mixed norms, and heat kernels}

We use the weak Lebesgue quasi-norm defined in the introduction. The interpolation estimate below is
a classical result. And it is the only weak-space interpolation which is needed in the proof.

\begin{lemma}[Weak-to-strong interpolation]\label{lem:interpolation}
If $1\le p<m<6$ and $v\in L^{p,\infty}\cap L^6$, then $v\in L^m$.
More precisely, if $1/m=\vartheta/p+(1-\vartheta)/6$, then
\begin{equation}\label{eq:interpolation}
 \|v\|_m\le C_{p,m}\|v\|_{L^{p,\infty}}^{\vartheta}
                    \|v\|_6^{1-\vartheta}.
\end{equation}
\end{lemma}

\begin{proof}
This is the Marcinkiewicz interpolation theorem applied to the identity map,
which is of weak type $(p,p)$ and of strong type $(6,6)$, see \cite{BerghLofstrom}.
\end{proof}

The order in every mixed norm considered in this paper is horizontal first and vertical second, that is
\begin{equation}\label{eq:mixed-norm}
 \|v\|_{L_z^\alpha L_h^c}
 =\left(\int_{\mathbb R}
       \left(\int_{\mathbb R^2}|v(x_h,z)|^c\,dx_h\right)^{\alpha/c}
       dz\right)^{1/\alpha}.
\end{equation}
The usual essential-supremum definition is used when an exponent is
infinite. 
The space $\dot H^1(\mathbb R^3)$ is the completion of
$C_c^\infty(\mathbb R^3)$ in the gradient $L^2$ norm.

\begin{lemma}\label{lem:mixed}
For $2<\alpha<\infty$ and $c=4\alpha/(\alpha-2)$, we have
\begin{equation}\label{eq:mixed-embedding}
 \|v\|_{L_z^\alpha L_h^c}\le C_\alpha\|\nabla v\|_2,
 \qquad v\in\dot H^1(\mathbb R^3).
\end{equation}
Consequently, if $6/5<m<2$ and $s=4m/(6-m)$, then
\begin{equation}\label{eq:mixed-product}
 \|va_h\|_{L_z^1L_h^s}
 \le C_m\|v\|_m\|\nabla a_h\|_2
\end{equation}
for $v\in L^m$ and $a_h\in\dot H^1(\mathbb R^3;\mathbb R^2)$.
\end{lemma}

\begin{proof}
Inequality \eqref{eq:mixed-embedding} is the mixed-norm Sobolev inequality of Benedek and Panzone \cite{BenedekPanzone}.
Estimate \eqref{eq:mixed-product} follows from it with $\alpha=m'$ and H\"older's inequality.
\end{proof}

We denote 
\begin{equation*}
    K_t^{(d)}(x)=(4\pi t)^{-d/2}\exp(-|x|^2/(4t))
\end{equation*}
to be the $d$-dimensional heat kernel, and write $K_t^h=K_t^{(2)}$.
For $t>0$, we usually write $e^{t\Delta}F=K_t^{(3)}*F$, and $e^{t\Delta_h}$, $e^{t\partial_z^2}$ are defined in an analogous way.
As we known, heat convolution opertors are bounded on $L^r$ for $1\le r\le\infty$.
Hence, the factorization $K_t^{(3)}=K_t^h\otimes K_t^{(1)}$ implies that
\begin{equation}\label{eq:heat-mixed}
 \|e^{t\Delta}F\|_{L_z^1L_h^s}
 \le\|F\|_{L_z^1L_h^s},\qquad 1\le s<\infty.
\end{equation}
That is, the operators are also bounded in the mixed Lebesgue space.

\subsection{Normalized pressure and the Hardy source}

We define the Riesz transform $R_j$ opertor by taking its Fourier symbol to be $-i\xi_j/|\xi|$. 
We use the classical $L^r$ boundedness of Riesz
transforms, the whole-space Calder\'on--Zygmund estimate for the Laplacian operator $\Delta$, and also the homogeneous Sobolev embedding, which are all the classical knowledge in harmonic analysis;
see \cite{Galdi,Stein1970}. What's more, We denote the real Hardy space by
$\mathcal H^1$, to distinguish it from the Sobolev space $H^1$.

\begin{lemma}[Background estimates]\label{lem:background}
Let $(u,b,P)$ be a smooth solution of \eqref{eq:mhd} satisfying
\eqref{eq:background}. Then $u,b\in\dot H^1\cap L^6\cap L^\infty$.
After subtracting a unique constant from the pressure,
\begin{equation}\label{eq:pressure-riesz}
 P=\sum_{i,j=1}^3R_iR_j(u_iu_j-b_ib_j)\in L^3,
 \qquad \nabla P\in L^{3/2}.
\end{equation}
Moreover,
\begin{equation}\label{eq:derivative-gain}
 \nabla^2u,\nabla^2b\in L^{3/2},
 \qquad \nabla u,\nabla b\in L^2\cap L^3.
\end{equation}
Finally, $S=-\Delta P$ belongs to $\mathcal H^1$ and
\begin{equation}\label{eq:hardy-source-bound}
 \|S\|_{\mathcal H^1}
 \le C\bigl(\|\nabla u\|_2^2+\|\nabla b\|_2^2\bigr).
\end{equation}
\end{lemma}

\begin{proof}
These are standard consequences of the whole-space Sobolev embedding, Calder\'on--Zygmund estimates, and the
Riesz-transform representation of the pressure.  The same analytic
package is used in the stationary MHD Liouville literature; see, for
example, \cite{Chae2023,ChenLiWang} and the references therein.
\end{proof}

Real Hardy space has many definition from different viewpoints. 
In this paper, We use the heat maximal characterization of real Hardy space
\cite{Stein1993}:
\begin{equation}\label{eq:hardy-maximal}
 \left\|\sup_{t>0}|e^{t\Delta}S|\right\|_1
 \le C\|S\|_{\mathcal H^1},
 \qquad \|S\|_1\le C\|S\|_{\mathcal H^1}.
\end{equation}

The Hardy--Littlewood maximal operator can be defined as
\begin{equation*}
\mathcal{M}f(x)=\sup_{r>0}\frac1{|B_r(x)|}\int_{B_r(x)}|f(y)|\,dy,
\end{equation*}
which is bounded on $L^r$ for $1<r<\infty$; heat averages satisfy
$|e^{t\Delta}g|\le C\mathcal{M}g$, and it holds uniformly in $t>0$.
These two maximal operators have distinct
roles: \eqref{eq:hardy-maximal}, rather than a strong $(1,1)$ bound
for $\mathcal{M}$, will control the pressure source.

% ===========================================================================
% Section 3 -- Structure lemma
% (formerly chapters/03_pressure_recovery.tex)
% ===========================================================================
\Needspace{10\baselineskip}
\section{Structure lemma}
\label{sec:pressure}
The lemma we will treat in this chapter is obtained by abstracting the estimates that we need in the process of proving Theorem~\ref{thm:main}. 
It is just a analytic lemma, independent of the
equations, which means proving this lemma itself will not require any information from MHD. 
Accordingly it is stated for the abstract data $f,F,E$ of \eqref{eq:directional}; in the MHD application they are
\[
f=\nu u_3,\qquad F=-(u_3u_h-b_3b_h),\qquad
E=\nu\partial_zu_3-u_3^2+b_3^2.
\]

\begin{proposition}[Finite pressure iteration]\label{prop:pressure}
Let $3/2<m<2$ and $s=4m/(6-m)$.
Suppose that
\begin{equation}\label{eq:pressure-hypotheses}
P\in L^3,\qquad S:=-\Delta P\in\mathcal H^1,\qquad
E\in L^2\cap L^3,\qquad f\in L^m,\qquad
F\in L_z^1L_h^s(\mathbb R^3;\mathbb R^2),
\end{equation}
and that, in distributions,
\begin{equation}\label{eq:directional}
\partial_z(P-E)=\Delta_h f+\operatorname{div}_h F.
\end{equation}
Then $P\in L^2$.
\end{proposition}

\begin{proof}
First, we use the action of the heat kernel operator to smooth the quantities and give them a high-low frequencies decomposition.
Put $Q=P-E$ and, for each of $P,E,Q,f,F$, denote convolution with the three-dimensional heat kernel $K_1^{(3)}$ by a subscript $0$.
Thus $Q_0=P_0-E_0$ and
\begin{equation}\label{eq:smoothed-directional}
\partial_zQ_0=\Delta_h f_0+\operatorname{div}_h F_0.
\end{equation}
Heat convolution is bounded on $L^m$ and also the mixed Lebesgue space by \eqref{eq:heat-mixed}.

Then, we treat the high frequency part simply.
Since the definition$-\Delta P=S$, we have
\[
\frac{d}{d\tau} e^{\tau \Delta}P = \Delta e^{\tau \Delta}P = -e^{\tau \Delta} S.
\]
Thus, intrgrating on $\tau \in [0,1]$, we obtain
\begin{equation}\label{eq:pressure-high-identity}
P-P_0=\int_0^1e^{\tau\Delta}S\,d\tau.
\end{equation}
Taking the $L^{\ell}$ norm on both sides. Using the Young's inequality, heat-kernel scaling and \eqref{eq:hardy-maximal}, we can compute that
\begin{equation}\label{eq:pressure-high-bound}
\|P-P_0\|_\ell
\le C\|S\|_1\int_0^1\tau^{-\frac32(1-1/\ell)}\,d\tau
\le C\|S\|_{\mathcal H^1}.
\end{equation}
So far, we do not confirm the exact range of the Lebesgue exponent $\ell$. 
However, our goal is to push the Lebesgue exponent of the pressure $P$ from 3 to 2, so $\ell$ must lie between them.
This fact guarantees the uniform boundedness of the constant $C$ above. The argument of high frequency part of $P$ has been finished.

Now, we are going to estimate the $P_0$ itself, which is more technical.
Let us first note that, if we endow $P$ with $L^r$ norm for some $2 \leq r \leq 3$, 
then interpolation gives
\begin{equation}\label{eq:Q-r-bound}
\|Q\|_r\le\|P\|_r+\|E\|_2+\|E\|_3,
\end{equation}
so $Q$ will be also in $L^r$.

We can observe that, the identity \eqref{eq:directional} controls the vertical derivative of $Q$, not $Q$ itself.
Therefore, our present task is to find some ways to obtain a pointwise estimate for it.
Now, we consider to apply to $Q$ a heat kernel operator carrying a parameter $\lambda$.
The role of this new heat operator is no longer high-low frequency decomposition as was done at the beginning of this proof. 
It now provides a device for absorbing derivatives since the quantities $f$ and $F$ just have integrable conditions, but have no condition about their derivatives. 
At the same time, it can supply a free parameter that may be chosen at will to sharpen our estimate.
This is something a mere $K_1^{(3)}$ convolution cannot achieve.
In detailed, for $0<\lambda\le1$, we set 
\[
w=e^{\lambda^{-2}\Delta_h}Q_0.
\]
For every interval $I=[z,z+T]$ of length parameter $T\ge 1$ and every $y\in I$, the Fundamental Theorem of Calculus gives that
\[
|w(x_h,z)|\le|w(x_h,y)|+\int_z^y|\partial_zw(x_h,s)|\,ds.
\]
This form is not suitable for an estimate, since the $w(x_h,z)$ on the left hand side and $w(x_h,y)$ on the right hand side carry almost the same information.
So we consider to take one more integral average over $I$ with respect to $y$, then we have
\begin{equation}\label{eq:vertical-average}
|w(x_h,z)|\le\frac1T\int_I|w(x_h,y)|\,dy
+\int_I|\partial_z w(x_h,y)|\,dy.
\end{equation}

Applying operator $e^{\lambda^{-2}\Delta_h}$ to \eqref{eq:smoothed-directional} gives
\begin{equation}\label{eq:vertical-derivative-solved}
\partial_z w
=e^{\lambda^{-2}\Delta_h}\partial_zQ_0
=\Delta_h\bigl(e^{\lambda^{-2}\Delta_h}f_0\bigr)
+\operatorname{div}_h\bigl(e^{\lambda^{-2}\Delta_h}F_0\bigr).
\end{equation}
Substituting \eqref{eq:vertical-derivative-solved} into \eqref{eq:vertical-average}, then we have explicitly
\begin{equation}\label{eq:three-term-splitting}
\begin{aligned}
|w(x_h,z)|\le{}&
\frac1T\int_I|w(x_h,y)|\,dy
+\int_I\bigl|\Delta_h\bigl(e^{\lambda^{-2}\Delta_h}f_0\bigr)(x_h,y)\bigr|\,dy\\
&+\int_I\bigl|\operatorname{div}_h\bigl(e^{\lambda^{-2}\Delta_h}F_0\bigr)(x_h,y)\bigr|\,dy .
\end{aligned}
\end{equation}
For the first and second terms of above (\ref{eq:three-term-splitting}), which kernels are $T^{-1}K^h_{\lambda^{-2}}\otimes\mathbf1_I$ and $\Delta_hK^h_{\lambda^{-2}}\otimes\mathbf1_I$ separately, by a scailing argument of Lebesgue space, we can straightly compute that
\begin{equation}\label{eq:average-kernel}
\big\|T^{-1}K^h_{\lambda^{-2}}\otimes\mathbf1_I\big\|_{r'}
=C_r\lambda^{2/r}T^{-1/r}
\end{equation}
and
\begin{equation}\label{eq:diffusion-kernel}
\big\|\Delta_hK^h_{\lambda^{-2}}\otimes\mathbf1_I\big\|_{m'}
=C_m\lambda^{2+2/m}T^{1-1/m}.
\end{equation}
For the last term $\operatorname{div}_h(e^{\lambda^{-2}\Delta_h}F_0)$, it can be written as the horizontal convolution of $F_0$ with $\nabla_hK^h_{\lambda^{-2}}$. Then H\"older's inequality gives
\begin{equation}\label{eq:flux-kernel}
\begin{aligned}
\int_I\big|\nabla_hK^h_{\lambda^{-2}}*_h F_0(x_h,y)\big|\,dy
&\le\|\nabla_hK^h_{\lambda^{-2}}\|_{s'}
\int_I\|F_0(\cdot,y)\|_{L_h^s}\,dy\\
&\le C_s\lambda^{1+2/s}\|F\|_{L_z^1L_h^s}.
\end{aligned}
\end{equation}
We can see there is no power of $T$ occurs in \eqref{eq:flux-kernel}. This means that using the mixed norm of $F$ has an effect on sharpening our estimate.
Combining \eqref{eq:three-term-splitting} with \eqref{eq:average-kernel}--\eqref{eq:flux-kernel}, and using $\|Q_0\|_r\le\|Q\|_r$ and $\|f_0\|_m\le\|f\|_m$, we obtain
\begin{equation}\label{eq:horizontal-bound}
\|e^{\lambda^{-2}\Delta_h}Q_0\|_\infty
\le C_{r,m}\big(\lambda^{2/r}T^{-1/r}\|Q\|_r
+\lambda^{2+2/m}T^{1-1/m}\|f\|_m
+\lambda^{1+2/s}\|F\|_{L_z^1L_h^s}\big).
\end{equation}

Since $\lambda\in(0,1]$, we want the power of $\lambda$ to be large. More precisely, we
want the three powers to share a common lower bound $\beta$.
The first term is $\lambda^{2/r}T^{-1/r}$, so the choice $T=\lambda^{-(\beta r-2)}$ makes its power exactly equals to $\beta$.
With this fixed $T$, the second power becomes $4-\beta r(1-1/m)$, which is at least $\beta$ exactly when $\beta\le4/(1+r(1-1/m))$.
And the last power is $1+2/s=1/2+3/m$, which exceeds $2$.
We therefore take
\begin{equation}\label{eq:beta-choice}
\beta=\min\left\{2,\frac4{1+r(1-1/m)}\right\},
\qquad T=\lambda^{-(\beta r-2)}.
\end{equation}
Indeed, since $r\le3$ and $m<2$, one has $8/5\le\beta\le2$ and $\beta r-2>0$, so $T\ge1$.
Therefore, with the new quantity
\[
X_r=\|P\|_r+\|E\|_2+\|E\|_3+\|f\|_m+\|F\|_{L_z^1L_h^s},
\]
the estimate becomes the following form
\begin{equation}\label{eq:horizontal-beta}
\|e^{\lambda^{-2}\Delta_h}Q_0\|_\infty\le C_{r,m}\lambda^\beta X_r.
\end{equation}
But $\lambda$ is still under a range restriction $0 < \lambda \leq 1$, whereas we want an estimate valid on the whole $(0,\infty)$.
For $t\ge1$, apply the $L^\infty$ boundedness of $e^{t\partial_z^2}$ to \eqref{eq:horizontal-beta} with $\lambda=t^{-1/2}$.
For $0<t<1$, use $Q_0=e^\Delta Q$ and the $L^r$--$L^\infty$ heat estimate.
It follows that
\begin{equation}\label{eq:negative-heat-bound}
B_\beta(Q_0):=\sup_{t>0}t^{\beta/2}\|e^{t\Delta}Q_0\|_\infty
\le C_{r,m}X_r.
\end{equation}
This is our final result on $w$.

The next step is where the Hardy space really comes into play. 
Let
\[
A(x)=\sup_{t>0}|e^{t\Delta}(-\Delta P_0)(x)|.
\]
Since $-\Delta P_0=e^\Delta S$, the Hardy maximal estimate (\ref{eq:hardy-maximal}) gives
\begin{equation}\label{eq:A-bound}
\|A\|_1\le C\|S\|_{\mathcal H^1}.
\end{equation}
The heat identity yields $|P_0-e^{t\Delta}P_0|\le tA$ holding in the sense of almost everywhere.
Writing $P_0=Q_0+E_0$, we obtain, for every $t>0$,
\begin{equation}\label{eq:pointwise-pressure}
|P_0(x)|\le tA(x)+B_\beta(Q_0)t^{-\beta/2}+C\mathcal{M}E_0(x).
\end{equation}
Optimizing in $t$ gives
\begin{equation}\label{eq:optimized-pressure}
|P_0(x)|\le C\mathcal{M}E_0(x)
+C_\beta A(x)^{\beta/(\beta+2)}B_\beta(Q_0)^{2/(\beta+2)}.
\end{equation}
Put
\begin{equation}\label{eq:improved-exponent}
a_m=\frac{m-1}{2m},\qquad
\ell=1+\frac2\beta
=\max\left\{2,\frac32+a_m r\right\}\in[2,9/4],
\end{equation}
so that $0<a_m<1/4$ because $1<m<2$.
Since we the only estimate on $A$ is related to $L^1$, we raise both sides of (\ref{eq:optimized-pressure}) to the power $\ell$ so that the powers match. This is the reason why we take $\ell$ in the form of (\ref{eq:improved-exponent}).
The $L^\ell$ boundedness of $\mathcal{M}$ and \eqref{eq:A-bound} imply
\begin{equation}\label{eq:pressure-low-bound}
\|P_0\|_\ell
\le C\|E\|_\ell
+C\|S\|_{\mathcal H^1}^{1/\ell}B_\beta(Q_0)^{1-1/\ell}.
\end{equation}

Thus, we obtain the starting point of the iteration, that is, together with \eqref{eq:pressure-high-bound} and \eqref{eq:negative-heat-bound}, this shows that $P\in L^r$ implies $P\in L^\ell$, with $\ell$ given by \eqref{eq:improved-exponent}.
The map $r\mapsto\ell$ of \eqref{eq:improved-exponent} is what we iterate.
Starting from $r_0=3$, which is available from Lemma~\ref{lem:background}, we set
\[
r_{n+1}=\max\left\{2,\frac32+a_m r_n\right\},
\]
by the improvement, $P\in L^{r_n}$ for every $n$.
The affine map $r\mapsto3/2+a_m r$ has slope $a_m$ and fixed point
\[
r_*=\frac{3/2}{1-a_m}=\frac{3m}{m+1}<2.
\]
Its untruncated iterates from $3$ are $r_*+a_m^n(3-r_*)$. Since $0<a_m<1$, they decrease geometrically to $r_*<2$, so after finitely many steps the untruncated value is below $2$. Therefore, we finally get $P \in L^2$ and finish the whole proof.
\end{proof}

% ===========================================================================
% Section 4 -- An L2 pressure criterion
% (formerly chapters/04_energy_criterion.tex)
% ===========================================================================
\section{\texorpdfstring{An $L^2$ pressure criterion}{An L2 pressure criterion}}\label{sec:energy}

In this chapter, we give a condition that directly implies $u=b=0$.
\begin{proposition}\label{prop:energy}
Let $(u,b,P)$ solve \eqref{eq:mhd} and satisfy \eqref{eq:background}.
If its normalized pressure belongs to $L^2$, then $u=b=0$.
\end{proposition}

\begin{proof}
Set $\zeta^\pm=u\pm b$, $\alpha=(\nu+\eta)/2$, and
$\gamma=(\nu-\eta)/2$. By adding and subtracting two equations in system \eqref{eq:mhd}, we get
\begin{equation}\label{eq:elsasser}
 -\alpha\Delta \zeta^\pm-\gamma\Delta \zeta^\mp
 +(\zeta^\mp\cdot\nabla)\zeta^\pm+\nabla P=0,
 \qquad \operatorname{div}\zeta^\pm=0.
\end{equation}

Then, we are going to find the test function to mimic the action of $u$ and $b$. For $\varepsilon>0$, define
\[
 T_\varepsilon(t)=\operatorname{sgn}(t)(|t|-\varepsilon)_+,
 \qquad
 \Psi_\varepsilon(t)=\tfrac12(|t|-\varepsilon)_+^2.
\]
The vanishing property of the function $(\cdot)_{+}$ implies that $T_\varepsilon(\zeta_i^\pm)$ and
$\Psi_\varepsilon(\zeta_i^\pm)$ have compact support for every fixed
$\varepsilon$. The chain rule gives
\begin{equation}\label{eq:truncation-gradient}
 \nabla T_\varepsilon(\zeta_i^\pm)
 =\mathbf{1}_{\{|\zeta_i^\pm|>\varepsilon\}}\nabla \zeta_i^\pm,
 \qquad
 \nabla\Psi_\varepsilon(\zeta_i^\pm)
 =T_\varepsilon(\zeta_i^\pm)\nabla \zeta_i^\pm.
\end{equation}
Consequently, $T_\varepsilon(\zeta_i^\pm)$ is an  compactly supported $H^1$ test function.
This suffices to verify that all the terms in \eqref{eq:elsasser} admits a well-defined
pairing with this test function by H\"older's inequality.

Then, the transport term vanishes exactly before taking any limit, that is
\[
 \int_{\mathbb R^3}(\zeta^\mp\cdot\nabla)\zeta_i^\pm
          T_\varepsilon(\zeta_i^\pm)\,dx
 =\int_{\mathbb R^3}\zeta^\mp\cdot\nabla\Psi_\varepsilon(\zeta_i^\pm)\,dx=0.
\]

Therefore, testing component $i$ of \eqref{eq:elsasser} now gives the following form
\begin{equation}\label{eq:truncated-energy}
\begin{aligned}
 0={}&\alpha\int_{\{|\zeta_i^\pm|>\varepsilon\}}|\nabla \zeta_i^\pm|^2\,dx
 +\gamma\int_{\{|\zeta_i^\pm|>\varepsilon\}}
             \nabla \zeta_i^\mp\cdot\nabla \zeta_i^\pm\,dx\\
 &-\int_{\{|\zeta_i^\pm|>\varepsilon\}}P\,\partial_i \zeta_i^\pm\,dx.
\end{aligned}
\end{equation}
A Sobolev function has zero gradient almost everywhere on each of its
level sets. This is a classc result, in our present smooth setting this follows directly
from the Implicit Function Theorem.
Thus, as $\varepsilon\to 0$, the gradient in \eqref{eq:truncation-gradient} converges to $\nabla \zeta_i^\pm$ strongly in $L^2$ space.
Products of two gradients are integrable, and
$P\,\partial_i \zeta_i^\pm\in L^1$ because $P\in L^2$.
Therefore, We may pass to the limit in all terms of \eqref{eq:truncated-energy}.

Sum first over $i$. The pressure terms cancel since
$\operatorname{div}\zeta^\pm=0$. Sum next over the two signs to obtain
\[
\begin{aligned}
 0&=\alpha\bigl(\|\nabla \zeta^+\|_2^2+\|\nabla \zeta^-\|_2^2\bigr)
       +2\gamma\int_{\mathbb R^3}\nabla \zeta^+:\nabla \zeta^-\,dx\\
  &=2\nu\|\nabla u\|_2^2+2\eta\|\nabla b\|_2^2.
\end{aligned}
\]
This implies $\nabla u = \nabla b = 0$ since $\nu, \eta$ are positive. Then, because of the Lebesgue integrability of $u$ and $b$, we obtain $u=b=0$.
\end{proof}

% ===========================================================================
% Section 5 -- Proofs of the main theorems
% (formerly chapters/05_proofs_main_theorems.tex)
% ===========================================================================
\section{Proofs of the main theorems}
\label{sec:main-proof}
\begin{proof}[Proof of Theorem~\ref{thm:main}]
Fix that
\begin{equation}\label{eq:m-choice}
\max\{p,q,3/2\}<m<2,\qquad s=\frac{4m}{6-m}\in(1,2).
\end{equation}
Lemma~\ref{lem:background} gives $u,b\in L^6$, then Lemma~\ref{lem:interpolation} gives
\begin{equation}\label{eq:component-integrability}
u_3,b_3\in L^m\cap L^4.
\end{equation}
Lemma~\ref{lem:background} shows the method to normalize the pressure $P$, and set
\begin{equation}\label{eq:directional-data}
f=\nu u_3,\qquad F=-(u_3u_h-b_3b_h),\qquad
E=\nu\partial_z u_3-u_3^2+b_3^2.
\end{equation}
Lemma~\ref{lem:structural} gives
\[
\partial_z P
=\nu\Delta u_3-\operatorname{div}_h(u_3u_h-b_3b_h)
-\partial_z(u_3^2-b_3^2).
\]

It suffices to verify the remaining hypotheses of Proposition~\ref{prop:pressure}.
First, $P\in L^3$ and $-\Delta P\in\mathcal H^1$ follow from Lemma~\ref{lem:background}.
Then, $f\in L^m$ by \eqref{eq:component-integrability}.
Third, the mixed product estimate \eqref{eq:mixed-product} straightly implies that
\begin{equation}\label{eq:main-flux-bound}
\|F\|_{L_z^1L_h^s}
\le C_m\big(\|u_3\|_m\|\nabla u_h\|_2
+\|b_3\|_m\|\nabla b_h\|_2\big)<\infty.
\end{equation}
Finally, $\partial_z u_3\in L^2\cap L^3$ by \eqref{eq:derivative-gain}, while $u_3^2,b_3^2\in L^2\cap L^3$ by their $L^4\cap L^6$ integrability.
These yield $E\in L^2\cap L^3$.
Proposition~\ref{prop:pressure} now yields $P\in L^2$, and Proposition~\ref{prop:energy} gives $u=b=0$. Thus $P$ is a constant by equations.
\end{proof}

\begin{proof}[Proof of Theorem~\ref{thm:elsasser-single}]
Without loss of generality, we consider the case $\zeta_3^+\in L^{p,\infty}$.
Repeatedly, fix
\[
\max\{p,3/2\}<m<2,
\qquad s=\frac{4m}{6-m}.
\]
Since $\zeta^\pm=u\pm b\in L^6$, Lemma~\ref{lem:interpolation} gives
$\zeta_3^+\in L^m\cap L^3\cap L^6$, while $\zeta^-\in L^6$.

For $\nu=\eta=\kappa$, the $\zeta^+$ equation is
\[
-\kappa\Delta \zeta^++(\zeta^-\cdot\nabla)\zeta^++\nabla P=0.
\]
The equations give the same identity in the Els\"asser variables:
\begin{equation}\label{eq:elsasser-directional-proof}
\partial_3(P-E)=\Delta_h f+\operatorname{div}_h F,
\qquad
f=\kappa \zeta_3^+,\qquad F=-\zeta_3^+\zeta_h^-,\qquad
E=\kappa\partial_3\zeta_3^+-\zeta_3^+\zeta_3^-.
\end{equation}
The mixed product estimate in Lemma~\ref{lem:mixed} yields
$F\in L_z^1L_h^s$.  Moreover,
\[
\|\zeta_3^+\zeta_3^-\|_2\le \|\zeta_3^+\|_3\|\zeta_3^-\|_6,
\qquad
\|\zeta_3^+\zeta_3^-\|_3\le \|\zeta_3^+\|_6\|\zeta_3^-\|_6,
\]
and the derivative estimate \eqref{eq:derivative-gain} gives
$\partial_3\zeta_3^+\in L^2\cap L^3$.  Thus $E\in L^2\cap L^3$.

This completes the verification of all the requisite conditions, likewise, we finish the proof.
\end{proof}

\begin{remark}
As what shows in the above proof, the reason why we only need weak $L^p$ property in one variable but not two in
the present setting is that, the term $F$ we have to treat has become simpler. This is nature, because the case $\nu = \eta$ is the simpler one compared with the general MHD.
\end{remark}

\begin{remark}\label{rem:endpoint}
The strict inequality $m<2$ in \eqref{eq:mixed-product} has a concrete reason:
the embedding fails at the endpoint.  For $R\ge e^2$ set
\[
 h_R(x_h)=
 \begin{cases}
  1,& |x_h|\le1,\\
  \log(R/|x_h|)/\log R,& 1<|x_h|<R,\\
  0,& |x_h|\ge R,
 \end{cases}
 \qquad
 v_R(x_h,z)=h_R(x_h)\chi(z/R),
\]
with $0\not\equiv\chi\in C_c^\infty(\mathbb R)$.  Direct calculation gives
$\|\nabla_hh_R\|_2^2=2\pi/\log R$ and $\|h_R\|_2^2\le CR^2/(\log R)^2$, hence
\[
 \|v_R\|_{L_z^2L_h^\infty}^2=R\|\chi\|_2^2,
 \qquad
 \|\nabla v_R\|_2^2
 =R\|\chi\|_2^2\|\nabla_hh_R\|_2^2
  +R^{-1}\|\chi'\|_2^2\|h_R\|_2^2
 \le\frac{CR}{\log R}.
\]
The quotient $\|v_R\|_{L_z^2L_h^\infty}/\|\nabla v_R\|_2$ is therefore at least
$c(\log R)^{1/2}$ and is unbounded, so no constant $C$ satisfies
$\|v\|_{L_z^2L_h^\infty}\le C\|\nabla v\|_2$ for all compactly supported
$v\in H^1(\mathbb R^3)$.
\end{remark}

% --------------------------- bibliography ----------------------------------

\end{document}